\documentclass[a4paper,12pt]{article}
\usepackage{graphicx}
\usepackage{mathptmx}
\usepackage{amssymb,amsmath,latexsym}
\usepackage[left=2cm,right=2cm, top=2cm,bottom=2cm,bindingoffset=0cm]{geometry}
\usepackage[T2A]{fontenc}
\usepackage{graphicx}
\usepackage{cite}
\usepackage{amssymb,amsmath,latexsym}
\usepackage{caption}
\usepackage{authblk}
\usepackage{color}

\usepackage{epsfig}

\begin{document}
\title{Synchronous oscillations and symmetry breaking in a model of two interacting ultrasound contrast agents}

\author[1]{Ivan R. Garashchuk}
\author[2]{ Alexey O. Kazakov}
\author[1]{Dmitry I. Sinelshchikov}
\affil[1]{National Research University Higher School of Economics, Moscow, Russia}
\affil[2]{National Research University Higher School of Economics, Nizhny Novgorod, Russia}

\maketitle

\begin{abstract}

We study nonlinear dynamics in a system of two coupled oscillators, describing the motion of two interacting microbubble contrast agents. In the case of identical bubbles, the corresponding symmetry of the governing system of equations leads to the possibility of existence of asymptotically stable synchronous oscillations. However, it may be difficult to create absolutely identical bubbles and, moreover, one can observe in experiments regimes that are unstable with respect to perturbations of equilibrium radii of bubbles. Therefore, we investigate the stability of various synchronous and asynchronous dynamical regimes with respect to the breaking of this symmetry. We show that the main factors determining stability or instability of a synchronous attractor are the presence/absence and the type of an asynchronous attractor coexisting with the synchronous attractor. On the other hand, asynchronous hyperchaotic attractors are stable with respect to the symmetry breaking in all the situations we have studied. Therefore, they are likely to be observed in physically realistic scenarios and can be beneficial for suitable applications when chaotic behavior is desirable.
\end{abstract}

\section{Introduction}
\label{intro}

In this work we study a nonlinear dynamical system consisting of two coupled forced nonlinear oscillators. This system describes the dynamics of two interacting microbubble contrast agents under the influence of an external periodic force. Microbubble contrast agents are micrometer size gas bubbles that are encapsulated into a visco-elastic shell. They are currently used for enhancing blood flow visualisation, see \cite{Szabo,Goldberg,Hoff} and there are also several possible further biomedical applications like noninvasive therapy and targeted drug delivery, see \cite{Klibanov2006,Coussios2008}.

It is known that nonlinear dynamics of microbubbles can be very complicated and depend on both control parameters and initial conditions (see, e.g. \cite{Parlitz1990,Behina2009,Macdonald2006,Carroll2013,Garashchuk2018,Garashchuk2019} and references therein). On the other hand, various types of bubbles dynamics can be both beneficial and undesirable depending on a particular application \cite{Hoff,Carroll2013}. Thus, it is important to thoroughly study the whole variety of possible stable dynamical regimes of contrast agents, transitions between them, and their dependence on the control parameters and initial conditions.

Here we concentrate on a model that describes oscillations of two microbubbles interacting via the Bjerknes force. This model was discussed in several works, see \cite{Takahira1995,Mettin1997,Ida2002,Doinikov2004,Dzaharudin2013}. Formally speaking, this model is a system of two coupled nonlinear oscillators with external periodic force. Some possible dynamical states and bifurcations occurring in this model were studied in \cite{Takahira1995,Macdonald2006,Dzaharudin2013}. A more detailed analysis of this system has been performed recently in \cite{Garashchuk2019}. In particular, it was shown that the dynamics of two coupled bubbles can be regular, quasiperiodic, chaotic and hyperchaotic and various multistable states are also possible. However, only symmetrical case, when both bubbles have the same equilibrium radius, was studied. On the other hand, the consideration of nonidentical bubbles may lead to the destruction of synchronous oscillations and to the emergence of new dynamical regimes. Investigation of the dynamics of nonidentical bubbles also seems more physically realistic since it would be technically difficult to create bubbles of exactly the same radius.

Moreover, as far as applications are concerned, if we analyze only synchronous oscillations while studying dynamics of a cluster of identical contrast agents, it is difficult to tell which of these regimes are physically realistic, because some of them may be unstable to the symmetry breaking. For example, in \cite{Behnia2019} authors study complicated dynamics in a cluster of interacting microbubbles with albumin and polymer shells. To simplify the calculations they use the GCM technique and study synchronous regimes only in a symmetrical situation (the identical bubbles case). The authors of \cite{Behnia2019} provide estimations of the stability domains of these regimes by certain parameters. However, it does not include the possibility of symmetry breaking in the system. If an attractor is unstable to small perturbations in bubbles' equilibrium radii, it means that it cannot appear in a system consisting of even slightly different bubbles. This means that such attractors are unlikely to be observed in a physically realistic situation. On the other hand, dynamical regimes that are stable to such perturbations are good candidates to be observed in a real ensemble of contrast agents.

Therefore, the main aim of this work is to study the influence of perturbations of bubbles' equilibrium radii on the dynamics in the model of two interacting gas bubbles. We show that there are three main possible scenarios of the impact of such perturbations. First, a multistable symmetrical state consisting of a synchronous and an asynchronous attractor can degenerate into a monostable asynchronous state. Second, the perturbations of equilibrium radii can lead to the appearance of new multistable states which do not exist in the symmetrical state. Third, some synchronous attractors can exist after perturbations of equilibrium radii leading to symmetry breaking. Further on we will call such attractors stable with respect to symmetry breaking or to the equilibrium radii perturbations.

It is important to clarify, that under stability of chaotic (hyperchaotic) attractors we do not mean their structural stability or even pseudohyperbolicity \cite{TS98,TS08,GGKK19,GKT20} --- a property that guarantee that every orbit in the attractor has one (two) positive Lyapunov exponents which never vanish under small perturbations of the system (parameter changing).  Throughout this work, under stability of the observed attractors we mean that these regimes under small perturbations in bubbles equilibrium radii are still observable in experiments. However, small stability windows inside them are not excluded.

We show that such stability becomes possible when a synchronous attractor is monostable or when a multistable state consisting of an asynchronous and a synchronous attractors is stable and does not degenerate into a monostable asynchronous state under symmetry breaking perturbations. What is more, we demonstrate that hyperchaotic attractors, which were recently observed in the studied model in \cite{Garashchuk2019}, are stable with respect to the symmetry breaking. Therefore, one can expect that such attractors will be experimentally observed. It also means that such attractors may be beneficial for applications where chaotic oscillations are desirable due to their stability.

The rest of this work is organized as follows. In Section \ref{sec:eqsys} we present the governing system of equations for the dynamics of two coupled microbubble contrast agents and introduce a non-dimensional parameter describing the symmetry breaking. In Section \ref{sec:one-route} we discuss various one-parametric routes of symmetry breaking and study the impact of multistability on the stability of synchronous attractors with respect to the symmetry breaking. In Section \ref{sec:2d} we perform two-dimensional analysis of the dynamics in the control parameters space and demonstrate that perturbations of bubbles equilibrium radii can lead to the emergence of new attractors. In the last Section we briefly summarize and discuss our results.

\section{Main system of equations}
\label{sec:eqsys}

\begin{figure}[!ht]
\center{
\includegraphics[width=0.5\linewidth]{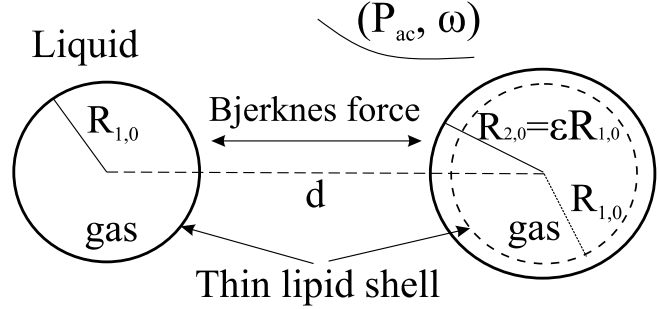}
}
\caption{Schematic picture of two interacting contrast agents oscillating in a liquid under the influence of an external pressure field. $\varepsilon$ denotes the ratio between equilibrium radii of two bubbles in static pressure of infinite volume of liquid without any external field, i.e. $R_{20} = \varepsilon R_{10}$, $|\varepsilon - 1| << 1$.}
\label{fig:UCA}
\end{figure}

We consider a model of oscillations of two interacting encapsulated gas bubbles in a liquid (see, Fig. \ref{fig:UCA}). If we denote time as $t$, bubbles radii as $R_{1}(t)$ and $R_{2}(t)$ respectively, derivatives with respect to $t$ as dots, the governing system of equations can be written in the following form:
\begin{eqnarray}
\label{eq:eq1}
    \left(1-\frac{\dot{R_{1}}}{c}\right)R_{1} \ddot{R_{1}} +\frac{3}{2}\left(1-\frac{\dot{R_{1}}}{3c}\right)\dot{R}_{1}^{2}=\frac{1}{\rho}\left[1+\frac{\dot{R}_{1}}{c}+\frac{R_{1}}{c}\frac{d}{dt}\right]P_{1}-\frac{d}{dt}\left(\frac{R_{2}^{2}\dot{R}_{2}}{d}\right), \vspace{0.1cm} \nonumber \\ 
     \left(1-\frac{\dot{R_{2}}}{c}\right)R_{2} \ddot{R_{2}} +\frac{3}{2}\left(1-\frac{\dot{R_{2}}}{3c}\right)\dot{R}_{2}^{2}=\frac{1}{\rho}\left[1+\frac{\dot{R}_{2}}{c}+\frac{R_{2}}{c}\frac{d}{dt}\right]P_{2}
     -\frac{d}{dt}\left(\frac{R_{1}^{2}\dot{R}_{1}}{d}\right), 
\end{eqnarray}
where
\begin{eqnarray*}
P_{i}=\left(P_{0}+\frac{2\sigma}{R_{i0}}\right)\left(\frac{R_{i0}}{R_{i}}\right)^{3\gamma}-\frac{4\eta_L \dot{R}_{i}}{R_{i}}-\frac{2\sigma}{R_{i}}-P_{0} - \quad\quad\quad\quad\quad\quad\quad\quad\quad\\
-4\chi\left(\frac{1}{R_{i0}}-\frac{1}{R_{i}}\right)-4\kappa_{S}\frac{\dot{R}_{i}}{R_{i}^{2}}-P_{ac}\sin(\omega t), \quad i=1,2.
\end{eqnarray*}
and $P_{\mathrm{stat}}$ is the static pressure, $P_v$ is the vapor pressure, $P_0 = P_{\mathrm{stat}} - P_v$, $P_{ac}$ is the magnitude of the external pressure field and $\omega$ is its cyclic frequency, $\sigma$ is the surface tension, $\rho$ is the density of the liquid, $\eta_{L}$ is the viscosity of the liquid, $c$ is the speed of sound in the liquid, $\gamma$ is the polytropic exponent, $\chi$ and $\kappa_{s}$ denote the shell elasticity and shell surface viscosity respectively.

Model \eqref{eq:eq1} consists of two generalized Raleigh--Plesset equations that are coupled via the Bjerknes forces, see Ref.-s \cite{Takahira1995,Mettin1997,Ida2002,Doinikov2004,Dzaharudin2013,Macdonald2006} for details. It also takes into consideration liquid's viscosity on the gas-liquid interface, liquid's compressibility via the Keller--Miksis model, see Ref. \cite{Keller1980}, surface tension and bubbles' shell according to the  de--Jong model, see Ref.-s \cite{deJong1992,Marmottant2005}.  The last term in the expression for $P_{i}$ describes the external pressure field that is assumed to be periodic.

Notice that system \eqref{eq:eq1} can be converted into an autonomous five-dimensional system of equations with respect to the following dependent variables $R_1, R_2, \dot R_1, \dot R_2$ and $\theta \in [0, 2 \pi]$. Further we use this five-dimensional system for numerical calculations, but we do not present it explicitly since it has a cumbersome form.

It can be seen that in the case of $R_{10} = R_{20}$ system \eqref{eq:eq1} possesses the following symmetry
\begin{equation}
R_1 \leftrightarrow R_2, \dot{R}_1 \leftrightarrow \dot{R}_2. \label{sym},
\end{equation}
which means that one can swap the indices between the bubbles. As a result, there exists an invariant manifold: $R_1(t) = R_2(t),\, \dot{R}_1(t) = \dot{R}_2(t)$. All synchronous solutions are embedded into this manifold and describe completely in-phase oscillations of the bubbles.  In Ref. \cite{Garashchuk2019} it was shown that solutions lying on this manifold can be both asymptotically attractive or repelling. Since the synchronization manifold is three-dimensional, synchronous asymptotically stable regimes can only be periodic (limit cycles) and simply chaotic \cite{Garashchuk2019} (with one positive Lyapunov exponent). On the other hand, attractors lying outside of the synchronization manifold can also be quasiperiodic and hyperchaotic (with two positive Lyapunov exponents) in addition to the previously mentioned types \cite{Garashchuk2019}. For asynchronous attractors in the system of two identical contrast agents, symmetry \eqref{sym} leads to two possibilities: either a given attractor is self-symmetrical with respect to the synchronization manifold or it has a symmetrical 'counterpart' lying on the 'other side' of the synchronization manifold.

In this work we are interested in scenarios of symmetry breaking via perturbation of the equilibrium radius, i.e. scenarios of destruction of the synchronization manifold. Without loss of generality, we consider perturbations of the equilibrium radius of the second bubble, while keeping the equilibrium radius of the first bubble constant. Consequently, we introduce the following non-dimensional parameter $\varepsilon$
\begin{equation}\label{eq:eq2}
  R_{20}=\varepsilon R_{10}, \quad  (\left| \varepsilon - 1 \right| \ll 1). 
\end{equation}
This type of perturbation is natural to a physically realistic situation, since it might be difficult to create an ensemble of exactly identical contrast agents. Furthermore, because the case of identical bubbles leads to some special properties, one can expect substantial qualitative changes in the behavior of system when the corresponding symmetry is lost. Synchronous attractors are primarily expected to undergo significant changes, since their properties are defined by the existence of the invariant manifold, which breaks down with an infinitely small perturbation in $R_{20}$. After such symmetry breaking two scenarios become possible: strange attractors which used to be synchronous can collapse. These attractors may also persist after such perturbations. We show in this paper that the evolution of atractors which used to belong to the invariant synchronization manifold significantly depends on presence and properties of the asynchronous attractors coexisting with a synchronous attractor.

Therefore, below we study the influence of small perturbations of the equilibrium radius of the second bubble on the nonlinear dynamics described by \eqref{eq:eq1}. We assume that $\varepsilon$, $P_{ac}$, $\omega$ and $d$ are the control parameters and the remaining parameters are fixed as follows: $P_v = 2.33$ kPa, $\sigma = 0.0725$ N/m, $\rho= 1000$ kg/m$^3$, $\eta_{L} = 0.001$ Ns/m$^3$, $c = 1500$ m/s, $\gamma = 4/3$, $\chi = 0.22$ N/m and $\kappa_S = 2.5\cdot 10^{-9}$ kg/s. These values of the parameters correspond to the adiabatic oscillations of two interacting SonoVue contrast agents with equilibrium radii $R_{i0}=1.72 \ \mu$m, $i=1,2$, see Ref. \cite{Tu2009}. To reduce the dimension of the control parameters space we also fix the frequency of the external pressure field at a physically relevant value of $\omega = 2.87 \cdot 10^7$ s$^{-1}$.

We perform all calculations in the following non-dimensional variables $R_{i}=R_{10}r_{i}$, $t=\omega_{0}^{-1}\tau$, where $\omega_{0}^{2}=3\kappa P_{0}/(\rho R_{10}^2)+2(3\kappa-1)\sigma/R_{10}+4\chi/R_{10}$ is the natural frequency of bubbles' oscillations. The non-dimensional bubbles speeds are given by $u_{i}=dr_{i}/d\tau=\dot{R}_{i}/ (R_{0}\omega_{0})$. We use the fourth-fifth order Runge--Kutta method for finding numerical solutions of the Cauchy problem for the considered system (see Ref. \cite{Cash1990} for details). For calculations of the Lyapunov spectra we use a standard algorithm (see Ref. \cite{Benettin1980}). Throughout this work the Lyapunov exponents are computed for the autonomous five-dimensional system. On all the graphs of Lyapunov exponents presented below we provide three exponents: $\lambda_1, \lambda_2, \lambda_f$, where $\lambda_1 \geq \lambda_2$ are the two effective exponents (they are responsible for the dynamics) and $\lambda_f\approx0$ is the referent exponent which is always close to zero and corresponds to the translations along an orbit. Computing colored two-parameter diagrams of Lyapunov exponents, it is important to separate effective Lyapunov exponents $\lambda_1$ and $\lambda_2$ from the slightly oscillating near zero referent exponent $\lambda_f$. For this purpose, we compare Lyapunov exponents not with zero, but with some small threshold $\lambda_{tr} = 0.001$. If $\lambda_1 < -\lambda_{tr}$, we conclude that the corresponding regime is periodic; if $-\lambda_{tr}<\lambda_1<\lambda_{tr}$, we mark this regime as quasiperiodic (however additional analysis can be necessary); if $\lambda_1 > \lambda_{tr}$ and $\lambda_2 < \lambda_{tr}$, we mark the regime as chaotic. Finally, if both exponents are greater than $\lambda_{tr}$ ($\lambda_i > \lambda_{tr},\, i=1,2$), we can guarantee that the corresponding regime is hyperchaotic. The Poincar\'e map is constructed by  taking values of phase variables at $t=k T, T=2\pi/\omega, k \in \mathbb{N}$, where $T$ is the period of the external pressure field.

While varying the control parameters, including $\varepsilon$, we also utilize the numerical continuation approach, inheriting the initial conditions step-by-step from certain starting values. If $\varepsilon$ is varied, we usually start from $\varepsilon = 1$ (if not specifically stated otherwise), taking initial conditions from attractors already found at a symmetrical case, and move in both directions by $\varepsilon$. At each step we skip transient processes which can be very long, especially after an attractor crisis. We compute the Lyapunov spectrum only after we believe that the trajectory has reached an attractor. If dynamics at $\varepsilon = 1$ is multistable,  we calculate the Lyapunov spectra for all the found coexisting attractors. More details on two-parameter Lyapunov analysis are given in Sec. 4.

\begin{figure}[!t]
\center{\includegraphics[width=1.0\linewidth]{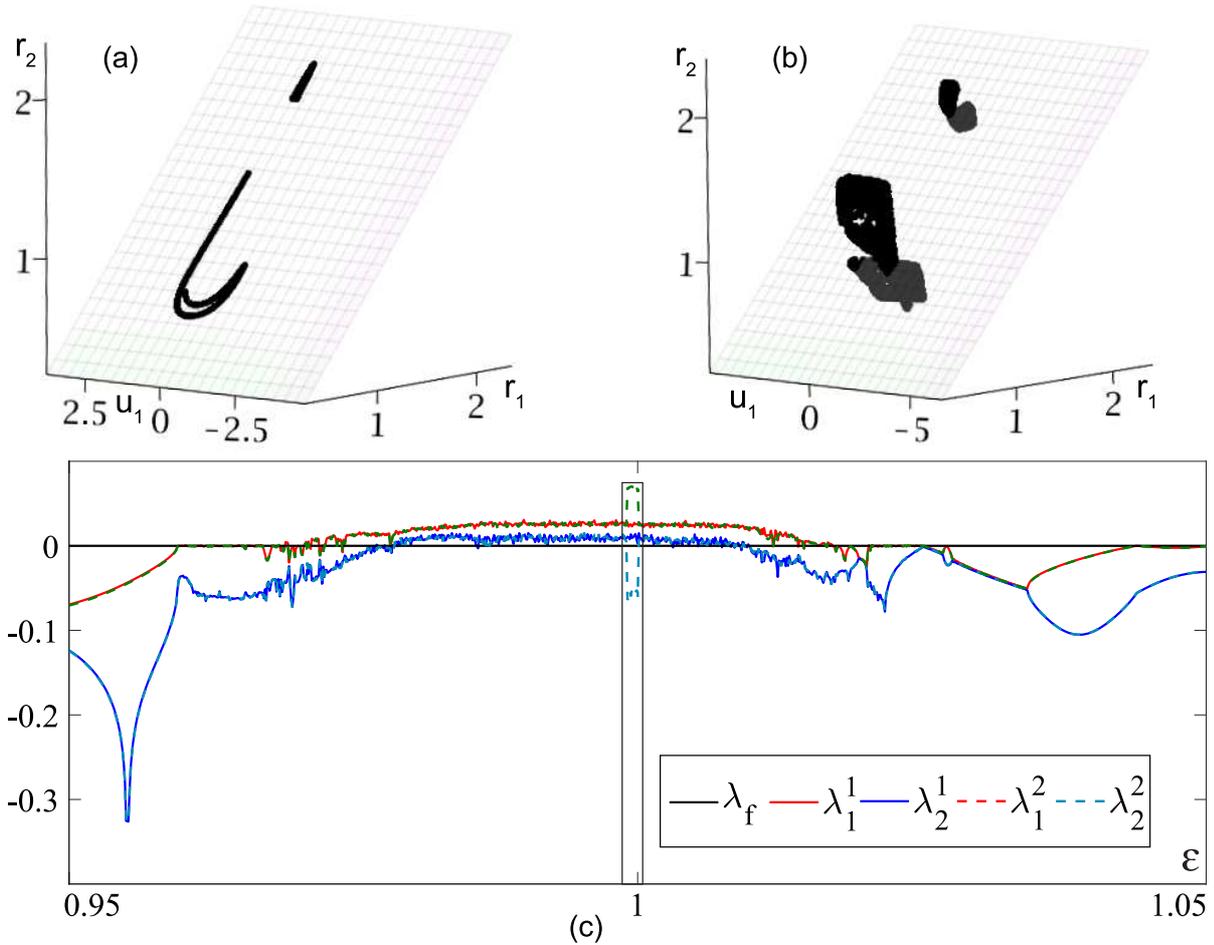}}
\caption{$P_{ac} = 1.2$ MPa, $d/R_{10} = 21$. 3-dimensional projection of Poincar\'e maps of (a) synchronous chaotic attractor, (b) asynchronous hyperchaotic attractor at $\varepsilon = 1$ on $r_1, u_1, r_2$, and (c) dependence of the Lyapunov spectra of attractors from (a) and (b) on $\varepsilon \in [0.95, 1.05]$}
\label{fig1}
\end{figure}

\section{One-parametric routes of symmetry breaking}
\label{sec:one-route}

In this section we study the influence of symmetry breaking on the bubbles' dynamics. We take $\varepsilon$ as a control parameter, fix the rest of the parameters and consider the impact of perturbations of bubble's equilibrium radius on the dynamics of the system.
Since usually different attractors can be efficiently distinguished by their spectra of Lyapnuov exponents, we calculate the dependence of the corresponding spectrum on $\varepsilon$ and demonstrate its three largest exponents. Notice also that some sharp transitions like an attractor crisis and jumps to another attractor in multistable states often can be seen from the dependence of the Lyapunov spectrum on the parameters.

\subsection{Symmetry breaking in multistable cases}
\label{ssec:multist}

\begin{figure}[!t]
\center{\includegraphics[width=\linewidth]{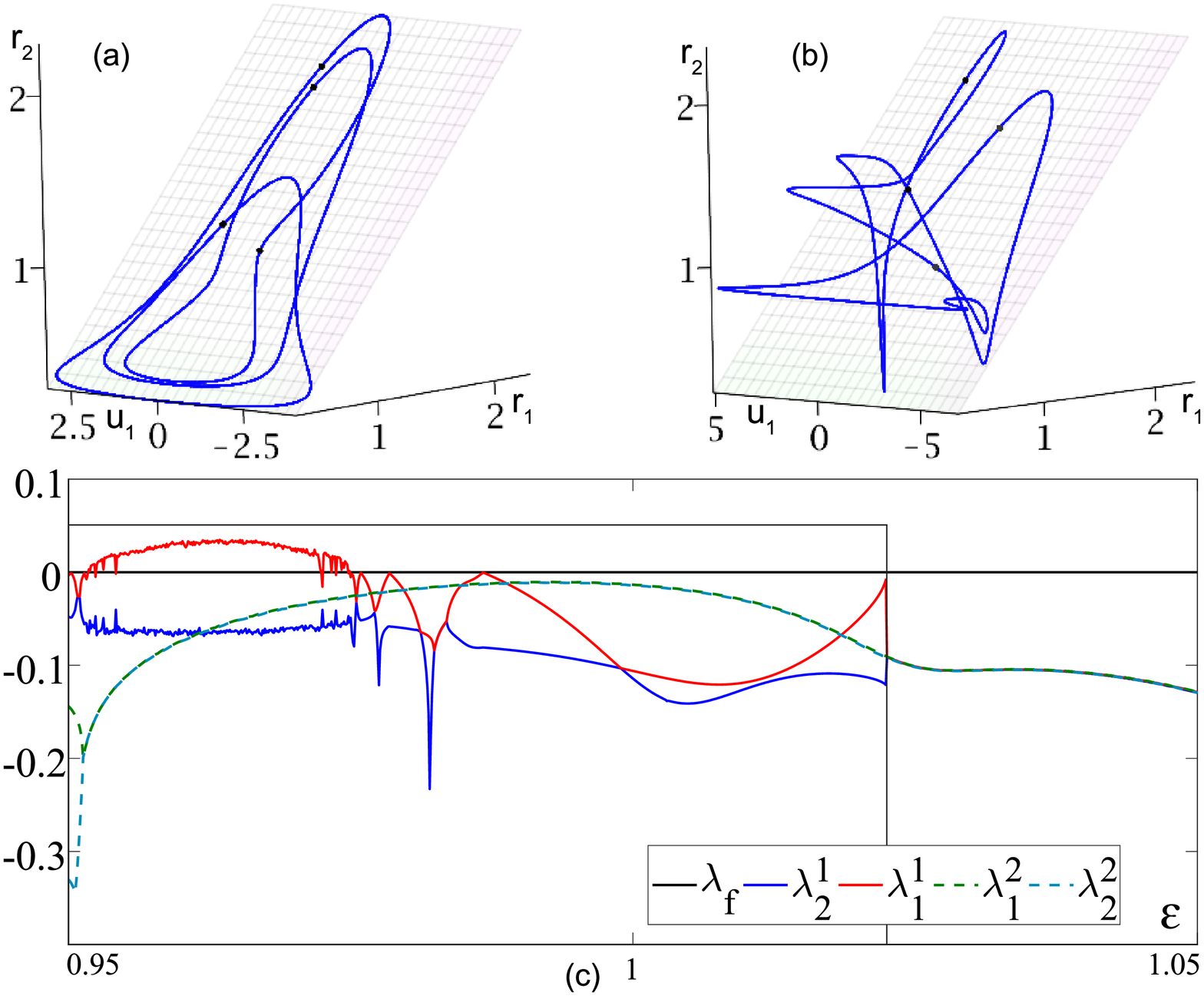}}
\caption{$P_{ac} = 1.2$ MPa, $d/R_{10} = 13$. 3-dimensional projections of the phase portraitrs (blue lines) and Poincar\'e maps (black dots) of (a) synchronous cycle of period 4, (b) asyncrhonous 4-periodic cycle at $\varepsilon = 1$ on $r_1, u_1, r_2$ , and (c)  dependence of the Lyapunov spectra of attractors (a) and (b) on $\varepsilon \in [0.95, 1.05]$.}
\label{fig2}
\end{figure}

Let us consider a multistable state discussed in \cite{Garashchuk2019} at $P_{ac} = 1.2$ MPa, $d/R_{10} = 21$, where two attractors coexist: a synchronous chaotic attractor and a hyperchaotic one (asynchronous), see Fig. \ref{fig1}a,b. Recall that a typical scenario of the emergence of a synchronous chaotic attractor in model \eqref{eq:eq1} is a cascade of period-doubling bifurcations of a limit cycle, while a hyperchaotic attractor typically appears according to the scenario presented in \cite{Garashchuk2019}.  Fig. \ref{fig1}c shows the dependence of two largest Lyapunov exponents on $\varepsilon$ for both attractors ($\lambda_i^1$, $i=1,2$ correspond to the hyperchaotic attractor (at $\varepsilon = 1$, see Fig. \ref{fig1}a) and $\lambda_i^2$, $i=1,2$ correspond to the synchronous chaotic attractor (at $\varepsilon = 1$, see Fig. \ref{fig1}b)). The strange attractor, which used to be synchronous at $\varepsilon = 1$, undergoes a crisis and orbits from its neighborhood to the asynchronous hyperchaotic attractor even if we subtly change $\varepsilon$, while the hyperchaotic attractor changes just slightly. If we consider larger regions of $\varepsilon$, the hyperchaotic attractor changes its type.

\begin{figure}[!h]
\center{\includegraphics[width=\linewidth]{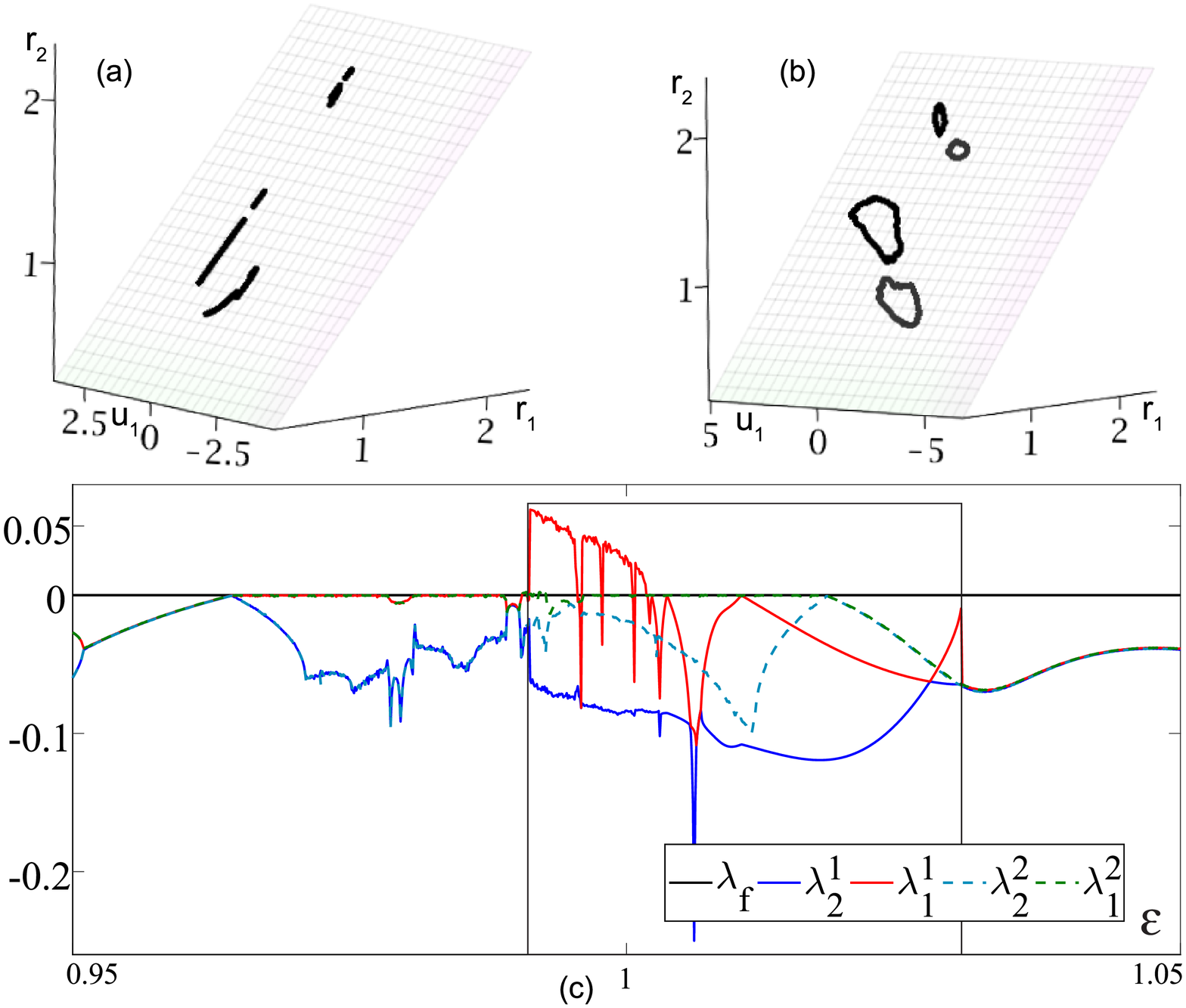}}
\caption{$P_{ac} = 1.2$ MPa, $d/R_{10} = 17.24$. 3-dimensional projection of the Poincar\'e map of (a) synchronous chaos, (b) torus at $\varepsilon = 1$ on $r_1, u_1, r_2$ and (c) dependence of the Lyapunov spectra of these attractors on $\varepsilon \in [0.95, 1.05]$.}
\label{fig3}
\end{figure}

On the other hand, a very different situation can be observed at the following values of the control parameters: $P_{ac} = 1.2$ MPa, $d/R_{10} = 13$. These values correspond to the coexistence of two limit cycles at $\varepsilon = 1$, one of which is synchronous and the other one is asynchronous, see Fig. \ref{fig2}a,b. In this case multistability is preserved for larger regions of $\varepsilon$, see Fig. \ref{fig2}c. The multistability window spreads till $\varepsilon = 0.95$ to the left, and to $\varepsilon \approx 1.0225$ to the right. This example demonstrates that in the case of the coexistence of an asynchronous limit cycle with a synchronous one, the later can be stable with respect to $\varepsilon$ (cf. the inverse situation of a synchronous chaotic attractor coexisting with a hyperchaotic one, Fig.-s \ref{fig1}c and \ref{fig2}c). The reason behind this is the type of an asynchronous attractor coexisting with synchronous one. Symmetry breaking allows coexistense of asynchronous limit cycles with the perturbed 'synchronous' limit cycles in quite a large range of $\varepsilon$.

To further illustrate this point, we consider the following values of the control parameters: $P_{ac} = 1.2$ MPa, $d/R_{10} = 17.24$. In the parameters space, these values lie between the values discussed in two previous cases. Now at $\varepsilon = 1$ we have attractors presented in Fig.\ref{fig3}a,b, which are obtained from those in Fig. \ref{fig2}a,b through the following bifurcation scenario. The synchronous 4-periodic cycle from Fig. \ref{fig2}a  goes through a period-doubling cascade and becomes the synchronous chaotic attractor demonstrated in Fig. \ref{fig3}a. The asynchronous cycle from Fig. \ref{fig2}b gives quasiperiodic regime after the Neimark--Sacker bifurcation, see Fig. \ref{fig3}b. The graphs of the Lyapunov spectra with respect to $\varepsilon$ for these attractors are shown in Fig. \ref{fig3}c. Here we observe a shorter range of values of $\varepsilon$ in which the attractor that used to be synchronous exists, than in the previous case (cf. Fig.-s \ref{fig3}c and \ref{fig2}c).

\begin{figure}[!h]
\center{\includegraphics[width=\linewidth]{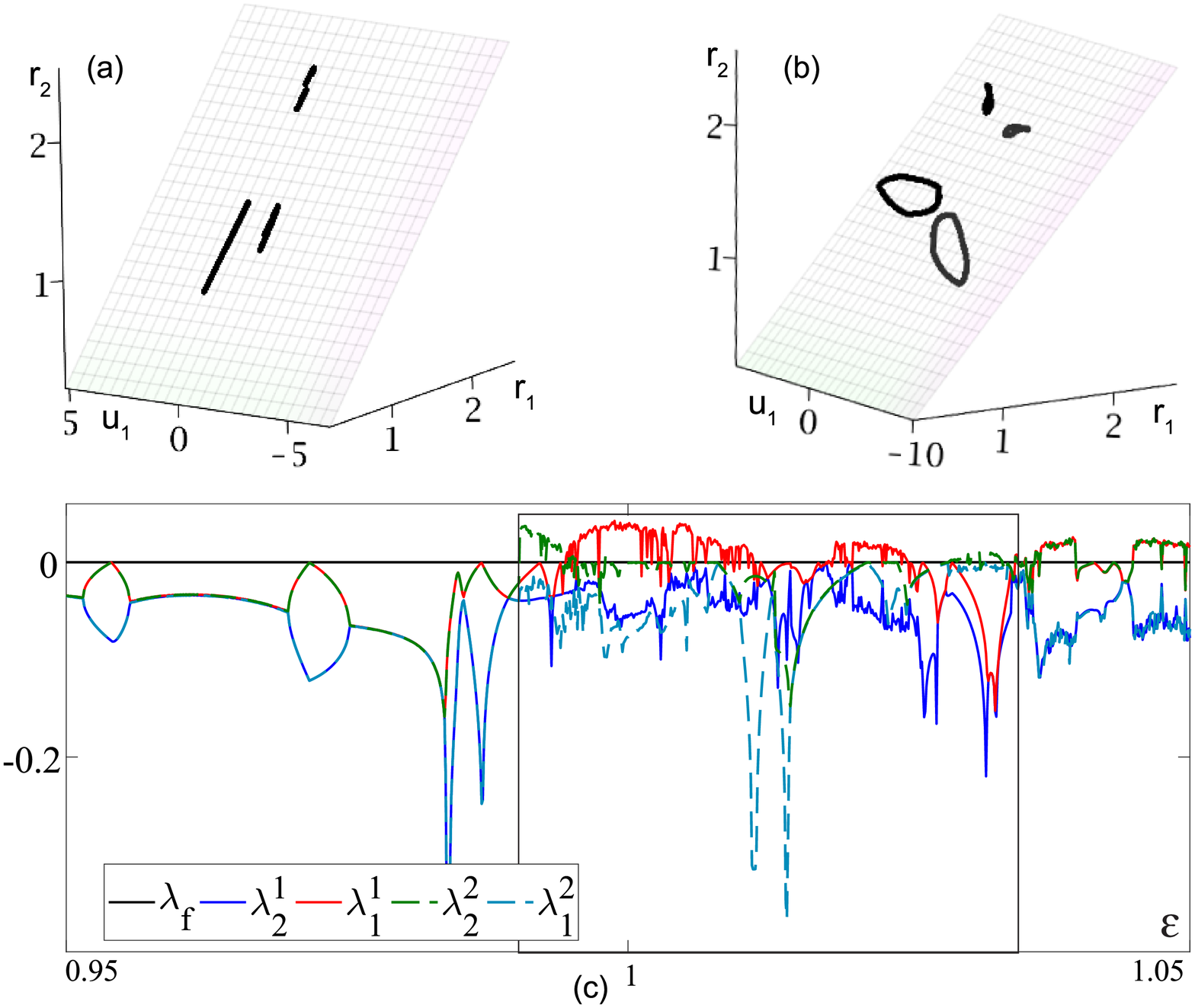}}

\caption{$P_{ac} = 1.2$ MPa, $d/R_{1,0} = 17.0$. 3-dimensional projection of the Poincar\'e map of (a) synchronous chaos, (b) torus on at  $\varepsilon = 1$ on $r_1, u_1, r_2$ and (c) dependence of the Lyapunov spectra of attractors (a) and (b) on $\varepsilon \in [0.95, 1.05]$.}
\label{fig4}
\end{figure}

A very similar situation can be observed for $P = 1.4$ MPa, $d/R_{10} = 17.0$, where a quasiperiodic regime coexists with a synchronous chaotic attractor like in the previous case, see Fig. \ref{fig4}a,b.
We demonstrate the dependence of the Lyapunov spectra on $\varepsilon$ in Fig. \ref{fig4}c. The area of multistability (stability of the synchronous attractor with respect to $\varepsilon$) looks similar to those in Fig. \ref{fig3}c. These results are consistent with our idea that the stability of a synchronous attractor under variations of symmetry breaking parameter $\varepsilon$ depends on the type of an asynchronous attractor.

\begin{figure}[!ht]
\center{\includegraphics[width=\linewidth]{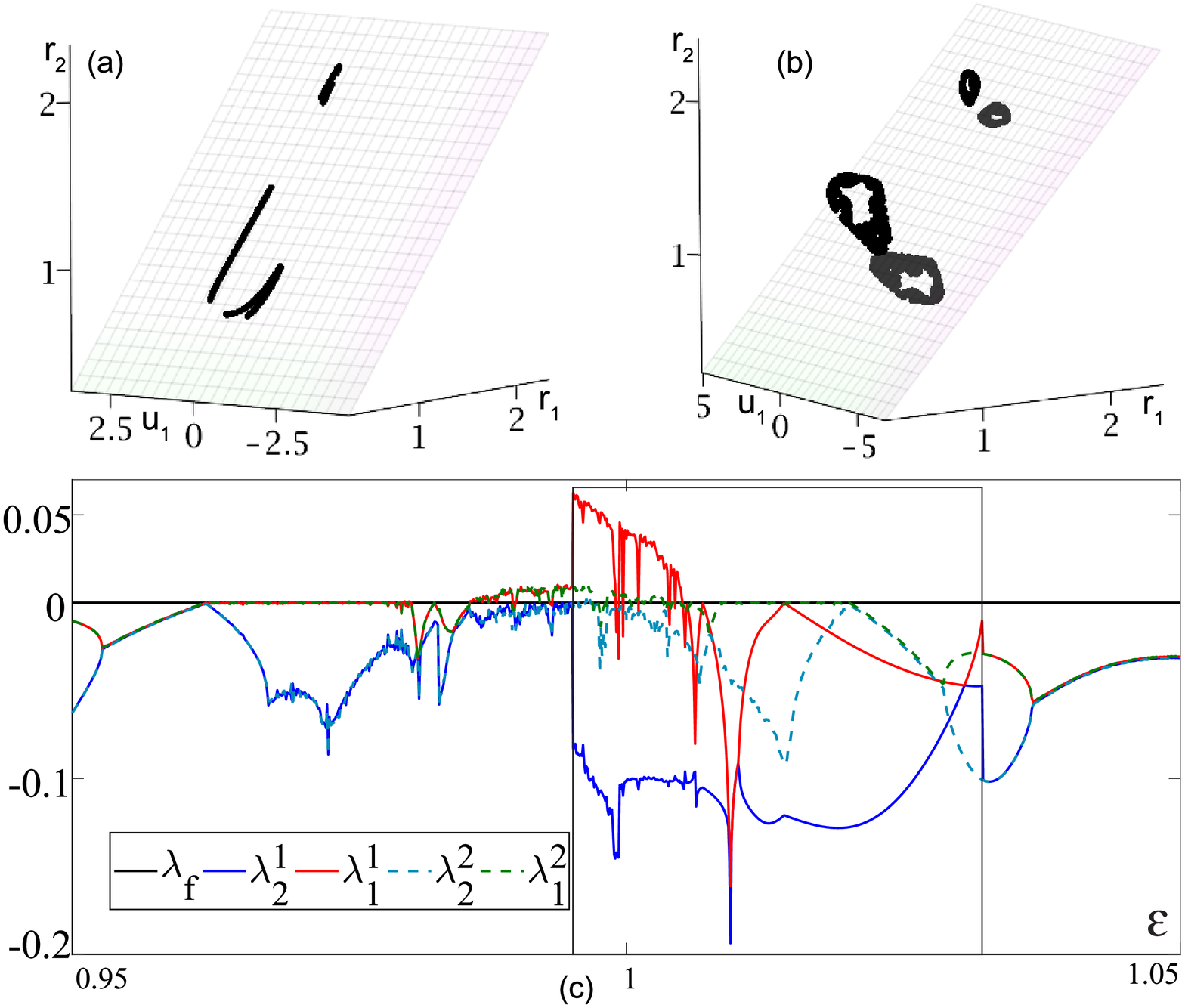}}
\caption{$P_{ac} = 1.2$ MPa, $d/R_{10} = 18.0924$. 3-dimensional projection of Poincar\'e map of (a) synchronous chaotic attractor, (b) asynchronous chaotic attractor at $\varepsilon = 1$ on $r_1, u_1, r_2$ and (c) dependence of the Lyapunov spectra of attractors (a), (b) on $\varepsilon \in [0.95, 1.05]$. }
\label{fig5}
\end{figure}

The next two types of attractors for which we perform similar analysis are a synchronous chaotic attractor and an asynchronous chaotic one at $P_{ac} = 1.2$ MPa, $d/R_{10} = 18.09$, see Fig. \ref{fig5}a,b. These attractors appear at the next stage in the bifurcation sequence happening to the attractors discussed earlier in Fig.-s \ref{fig2}a,b and \ref{fig3}a,b. The asynchronous attractor is now chaotic: stable limit cycle inside a resonant region on a torus presented in Fig. \ref{fig5}b undergoes the secondary Neimark--Sacker bifurcation, then secondary torus evolves into a torus-chaos attractor via a cascade of period doubling bifurcation with high-periodic resonant orbit.
The plots of the Lyapunov spectra of these attractors with respect to $\varepsilon$ are presented in Fig. \ref{fig5}c. The range of stability of synchronous chaotic attractor is very similar to those demonstrated in Fig. \ref{fig3}c. Therefore, we see that the impact of a torus-chaos attractor on the stability with respect to $\varepsilon$ of a synchronous chaotic attractor is similar to that of a quasiperiodic attractor.

\begin{figure}[!ht]
\center{\includegraphics[width=\linewidth]{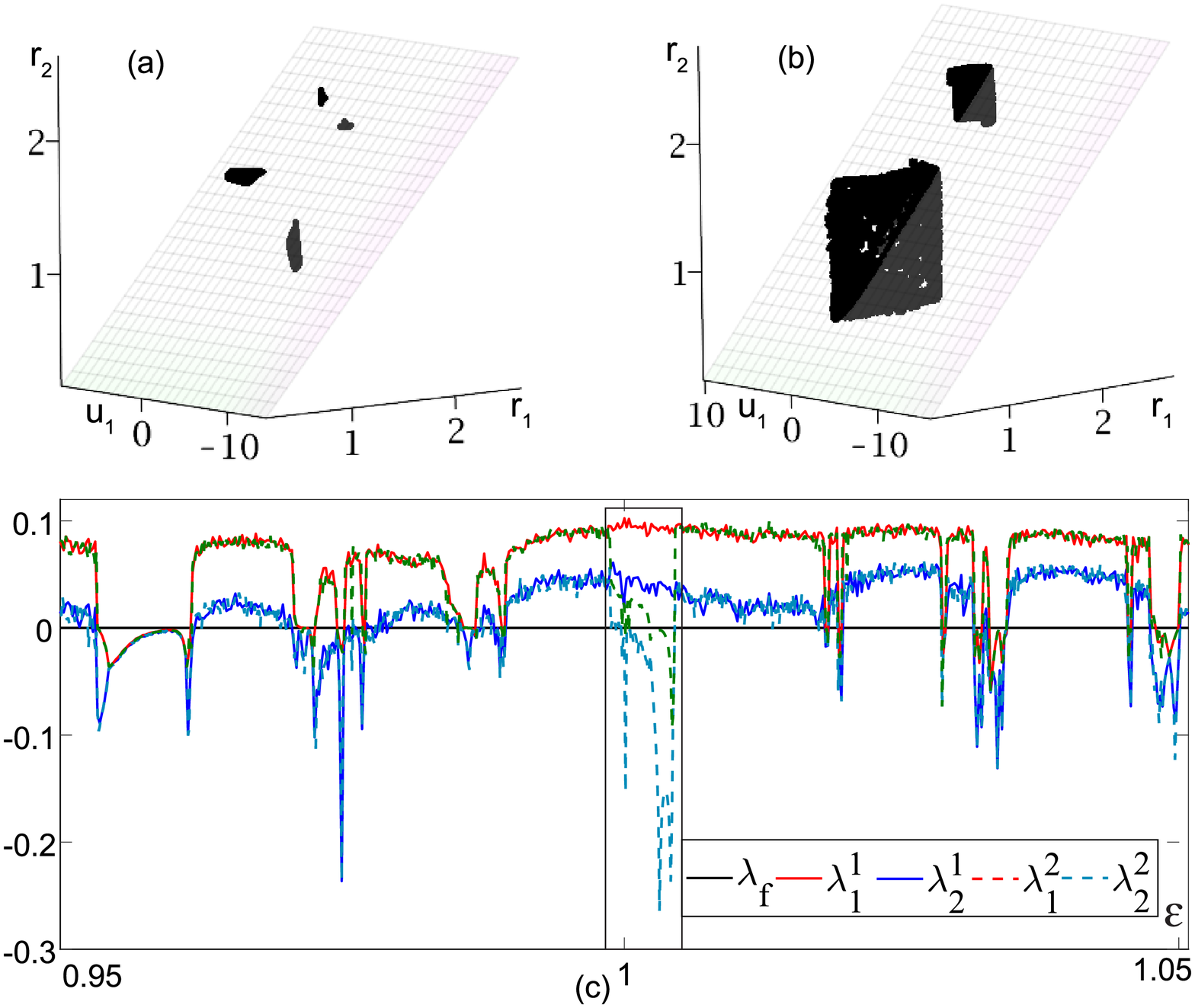}}
\caption{ $P_{ac} = 1.68$ MPa, $d/R_{10} = 24.93$. 3-dimensional projection of Poincar\'e map  of (a) synchronous chaotic attractor, (b) hyperchaotic attractor at  $\varepsilon = 1$ on $r_1, u_1, r_2$ and (c) dependence of the Lyapunov spectra of attractors (a), (b) on $\varepsilon \in [0.95, 1.05]$. }
\label{fig6}
\end{figure}

Now let us consider another type of multistable state at $P_{ac} = 1.68$ MPa, $d/R_{10} = 24.93$, where coexistence of asynchronous chaotic and hyperchotic attractors takes place, see Fig. \ref{fig6}a,b. This asynchronous chaotic attractor emerges via the Afrai\-movich--Shilnikov scenario, see Ref. \cite{AfrShil1983}. In Fig. \ref{fig6}c we can observe that there is a very narrow range of $\varepsilon$ where these attractors coexist. The asynchronous chaotic attractor disappears after very small perturbations of $\varepsilon$ and outside of this interval in $\varepsilon$ only the hyperchaotic attractor exists. This example once again demonstrates the stability of hyperchaotic attractors in multistable states. On the other hand, if an attractor coexists with a hyperchaotic one at $\varepsilon = 1$, it is highly likely that it will disappear with small perturbations of $\varepsilon$. This particular example shows that even some asynchronous attractors can be unstable with respect to the symmetry breaking if they coexist with hyperchaotic regimes.

\begin{figure}[!ht]
\center{\includegraphics[width=\linewidth]{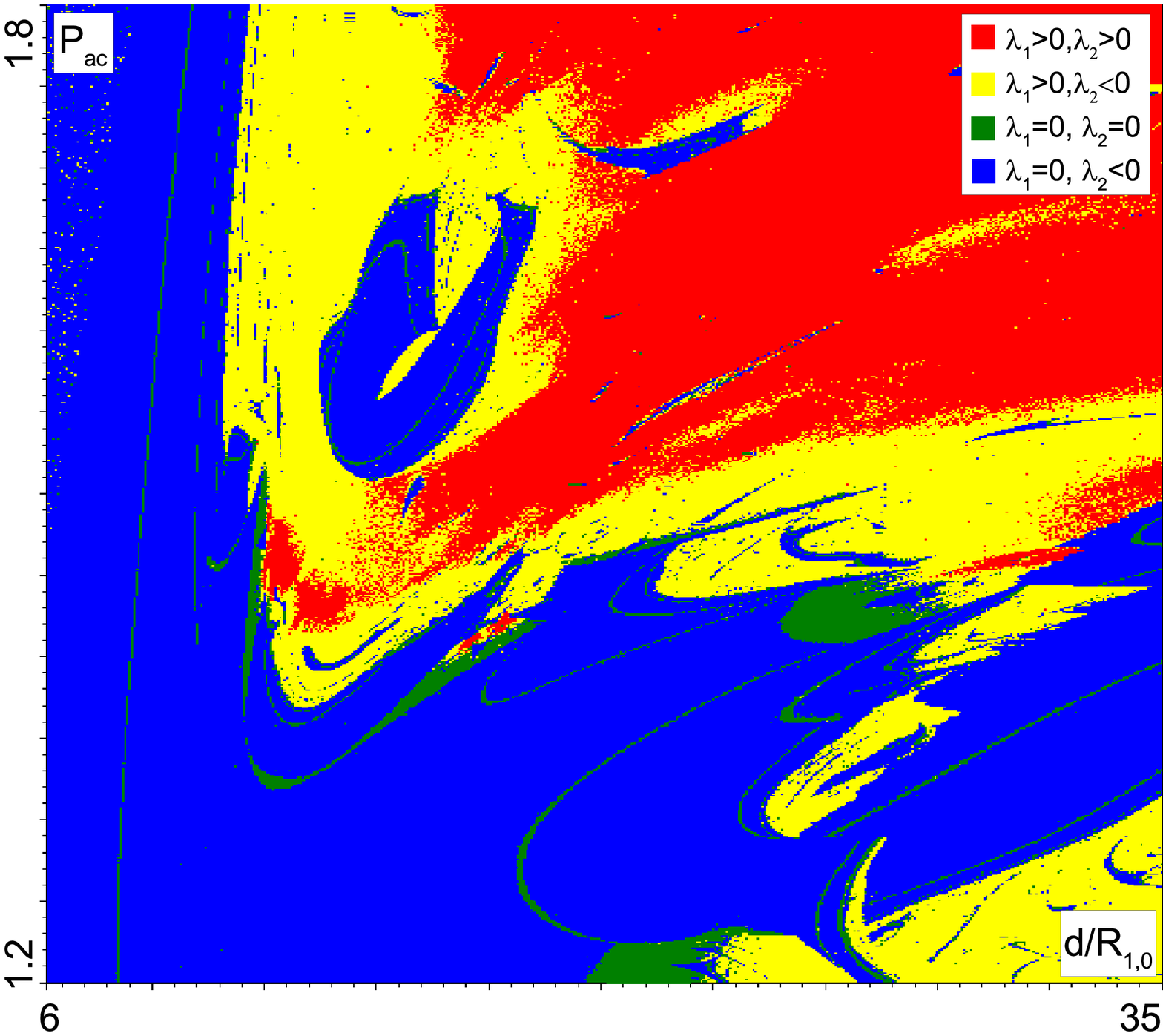} }
\caption{Chart of the Lyapunov exponents for the following parameters area: $P_{ac} \in [1.2, 1.8]$ MPa, $d/R_{10} \in [6, 35]$, $\varepsilon = 1.024$. }
\label{fig:2d_map_1}
\end{figure}

\section{Two-dimensional analysis}
\label{sec:2d}

Stability of hyperchaotic attractors is also demonstrated in two-dimensional chart of dynamical regimes in asymmetrical case: $\varepsilon = 1.024$, see Fig. \ref{fig:2d_map_1}. For computing this chart we start from the point $(r_1, u_1, r_2, u_2) =  (1.09, -0.47, 0.77, 0.49)$ corresponding to the asynchronous hyperchaotic attractor at $(d/R_{10}, P_{ac}) = (17.5, 1.52 \times 10^6)$. First, we continue this regime up and down to the interval $P_{ac} \in [1.2  \times 10^6, 1.8  \times 10^6]$, and, then, to the left and to the right varying $d/R_{10} \in [6,35]$.

In the symmetrical case multistable states used to take place in a lot of areas in the control parameters space. Now hyperchaos occupies a lot of those areas, for which one of the coexisting attractors was hyperchaotic. Thus we can conclude that hyperchaotic oscillations tend to be stable with respect to symmetry breaking. Moreover, they are likely to be observed in physical systems with the corresponding values of the control parameters. It means that if chaotic behavior of microbubbles is desirable, it is reasonable to use control parameters from the areas in which a hyperchaotic attractor exists.

\begin{figure}[!ht]
\center{
\includegraphics[width=\linewidth]{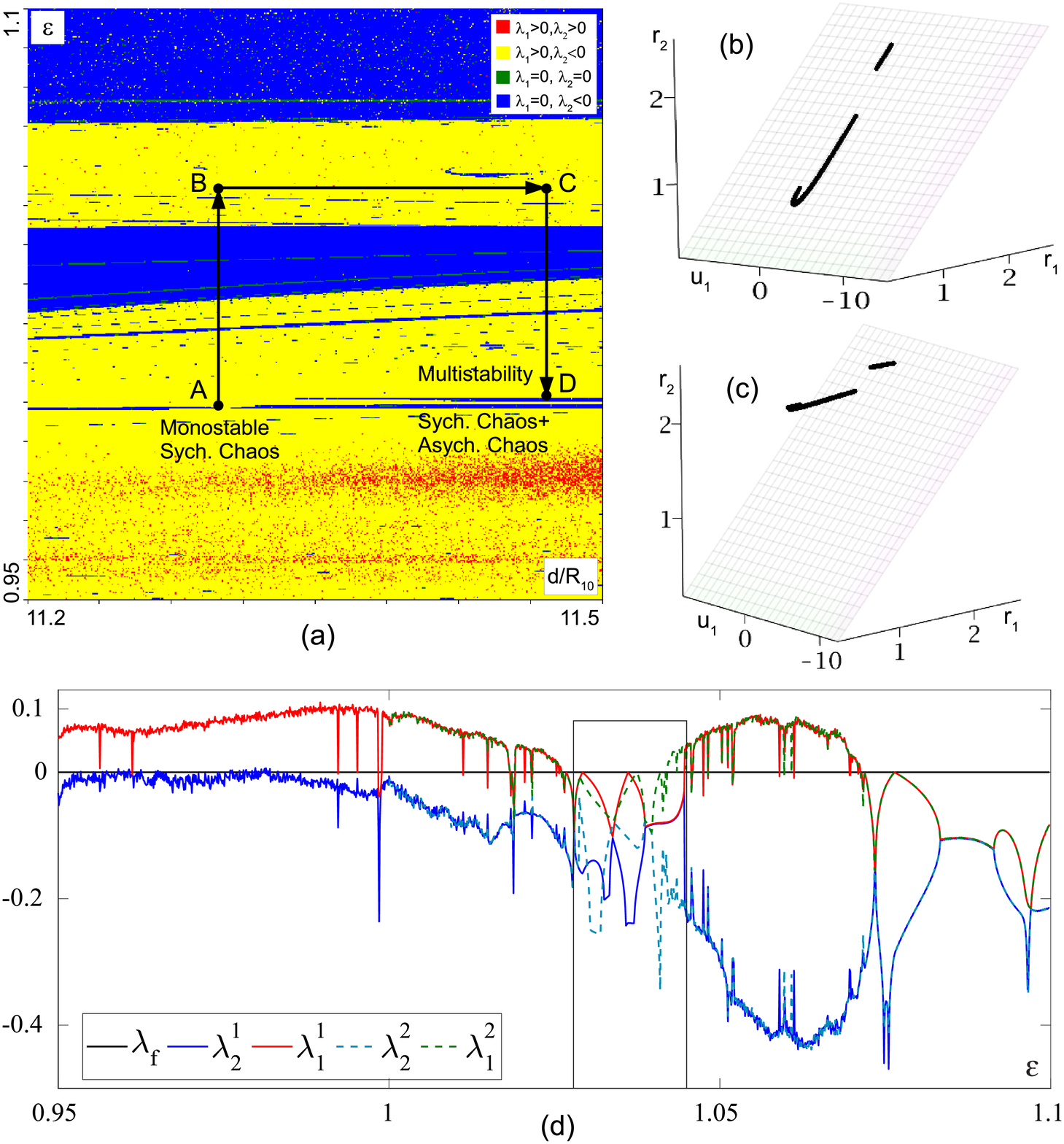}
}

\caption{$P_{ac} = 1.8$ MPa. (a)chart of Lyapunov exponents in the parameters area $d/R_{10} \in [11.2, 11.5]$, $\varepsilon \in [0.95, 1.1]$, (b) and (c) 3-d projections of the Poincar\'e section of (b) monostable synchronous chaotic attractor at $ d/R_{10} = 11.3, \varepsilon = 1$ (point A) and (c) new asymmetrical chaotic attractor at $d/R_{10} = 11.3, \varepsilon = 1.05$ (point B) on $r_1, u_1, r_2$, (d) graph of the Lyapunov spectra of attractors (b), (c) at $d/R_{10} = 11.3$ by $\varepsilon \in [0.95, 1.1]$ (route AB). }

\label{fig8}
\end{figure}

\subsection{Switch of synchronous attractors' stability}
\label{sec:switch}

Let us consider the following 2-dimensional area:  $P_{ac} = 1.8$ MPa, $d/R_{10} \in [11.2, 11.5]$, $\varepsilon \in [0.95, 1.1]$. We demonstrate a 2-dimensional chart of Lyapunov exponents of one of stable branches of dynamical regimes in Fig. \ref{fig8}a. For computing this diagram we start from the point $(r_1, u_1, r_2, u_2) =  (0.924, 1.118, 0.924, 1.118)$ corresponding to the synchronous chaotic attractor at $(d/R_{10}, \varepsilon) = (11.3, 1)$. We continue this attractor up and down to the interval $\varepsilon \in [0.95, 1.1]$, and, then, to the left and to the right varying $d/R_{10} \in [11.3, 11.5]$.

After all the previous analysis of symmetry breaking in multistable states, it is quite interesting to study the influence of a symmetry breaking on a monostable synchronous attractor. For instance, such a state occurs at point A in Fig. \ref{fig8}a: $P_{ac} = 1.8$ MPa, $d/R_{10} = 11.3$, $\varepsilon = 1$. At these values of parameters a monostable synchronous attractor exists, see Fig. \ref{fig8}b. This synchronous attractor is stable with respect to perturbations in $\varepsilon$. It changes smoothly with small perturbations in $\varepsilon$ and exists in quite a large range of $\varepsilon$, see Fig. \ref{fig8}d (route AB in Fig. \ref{fig8}a). Note that while for small $\varepsilon$ this attractor is monostable, at higher values of $\varepsilon$ a new attractor emerges, which exists only in asymmetrical states at relatively high values of $\varepsilon$, see Fig. \ref{fig8}c. The multistability window in Fig. \ref{fig8}d corresponds to the coexistence of this new attractor with previously synchronous one. If $\varepsilon$ is further increased, the previously synchronous attractor disappears, and new asymmetrical one becomes the only remaining attractor that we could find.

\begin{figure}[!ht]
\center{\includegraphics[width=\linewidth]{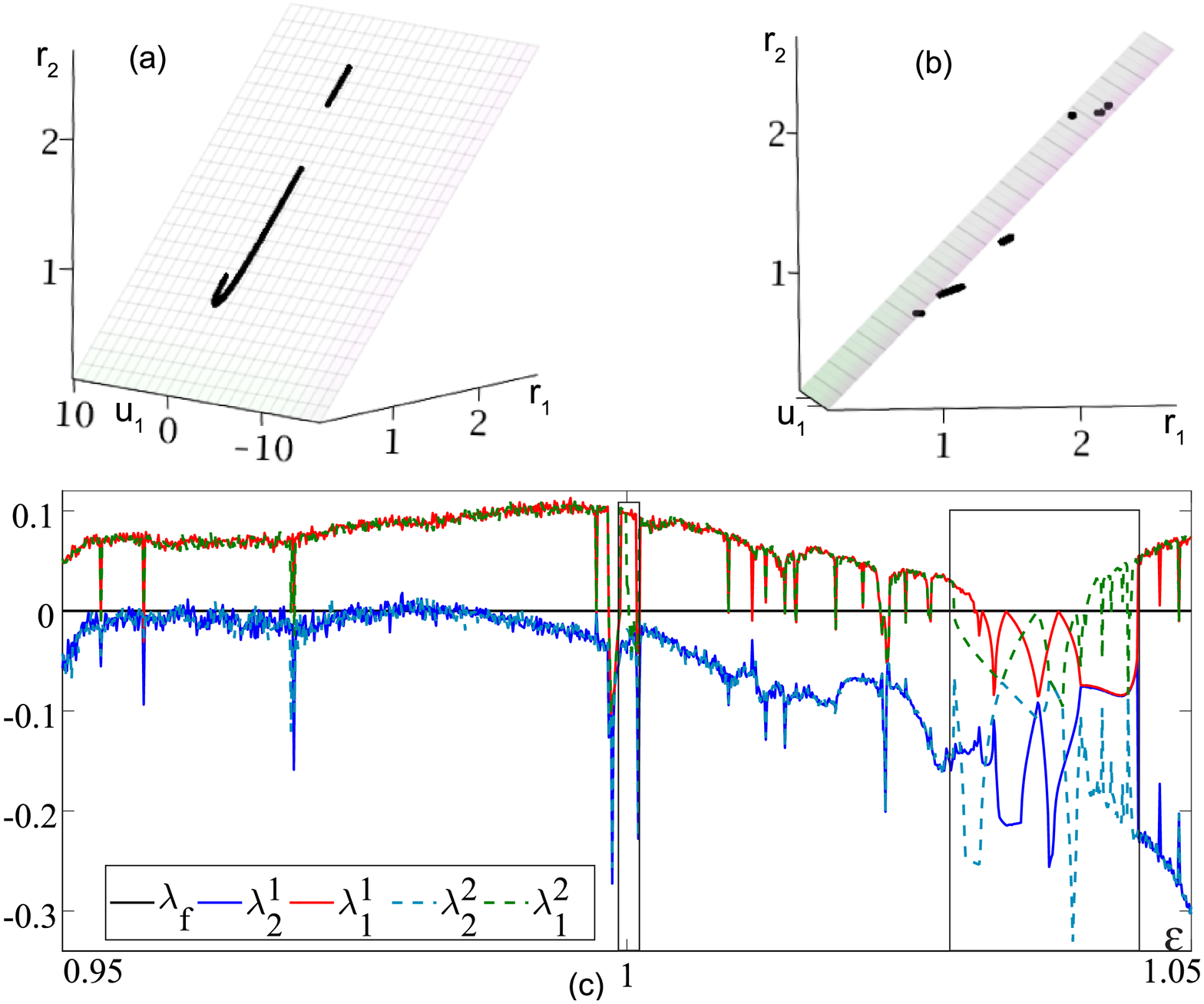}}

\caption{$P_{ac} = 1.8$ MPa, $d/R_{10} = 11.47$. 3-d projections of the Poincar\'e section of (a) synchronous chaotic attractor, (b) asynchronous chaotic attractor at $ d/R_{10} = 11.47, \varepsilon = 1$ (point D in Fig. \ref{fig8}a) on $r_1, u_1, r_2$, (c) graph of dependence of the Lyapunov spectra of attractors at $d/R_{10} = 11.47$ on $\varepsilon \in [0.95, 1.05]$ (route CD in Fig. \ref{fig8}a). }
\label{fig9}
\end{figure}

On the other hand, at $P_{ac} = 1.8$ MPa, $d/R_{10} = 11.47$, $\varepsilon = 1$ (point D in Fig. \ref{fig8}a), there coexist two attractors: synchronous chaotic and asynchronous chaotic, see Fig. \ref{fig9}a,b. One can notice that the synchronous attractor in Fig. \ref{fig9}a looks almost exactly the same as the one in Fig. \ref{fig8}b, which is obtained with the help of the numerical continuation method by varying the parameter $d/R_{10}$.

Let us apply numerical continuation for both attractors at $P_{ac} = 1.8$ MPa, $d/R_{10} = 11.47$, $\varepsilon = 1$ (point D in Fig. \ref{fig8}a) and for the new asymmetrical attractor at $\varepsilon = 1.05$ (point C in Fig. \ref{fig8}a). As a result, varying $\varepsilon$, we can find the dependence of the Lyapunov spectra on $\varepsilon$ for all the attractors existing for these parameters in the interval $0.95 < \varepsilon < 1.05$, see Fig. \ref{fig9}c. One can observe two windows of multistability: first very narrow one around $\varepsilon = 1$, and the second one at larger values of $\varepsilon$. The first one corresponds to the coexistence of synchronous and asynchronous chaotic attractors at $\varepsilon = 1$. Since this interval is very narrow, we can conclude that the synchronous attractor is very sensitive to variations in $\varepsilon$, and, consequently, is not physically realistic. The second multistability window corresponds to the coexistence of previously asynchronous attractor with the new asymmetrical one. The comparison of two attractors in the second multistability window also demonstrates that the new asymmetrical attractor is not a continued form of the asynchronous attractor existing at $P_{ac} = 1.8$ MPa, $d/R_{10} = 11.47$, $\varepsilon = 1$ (point D in Fig. \ref{fig8}a).

\begin{figure}[!ht]
\center{\includegraphics[width=0.9\linewidth]{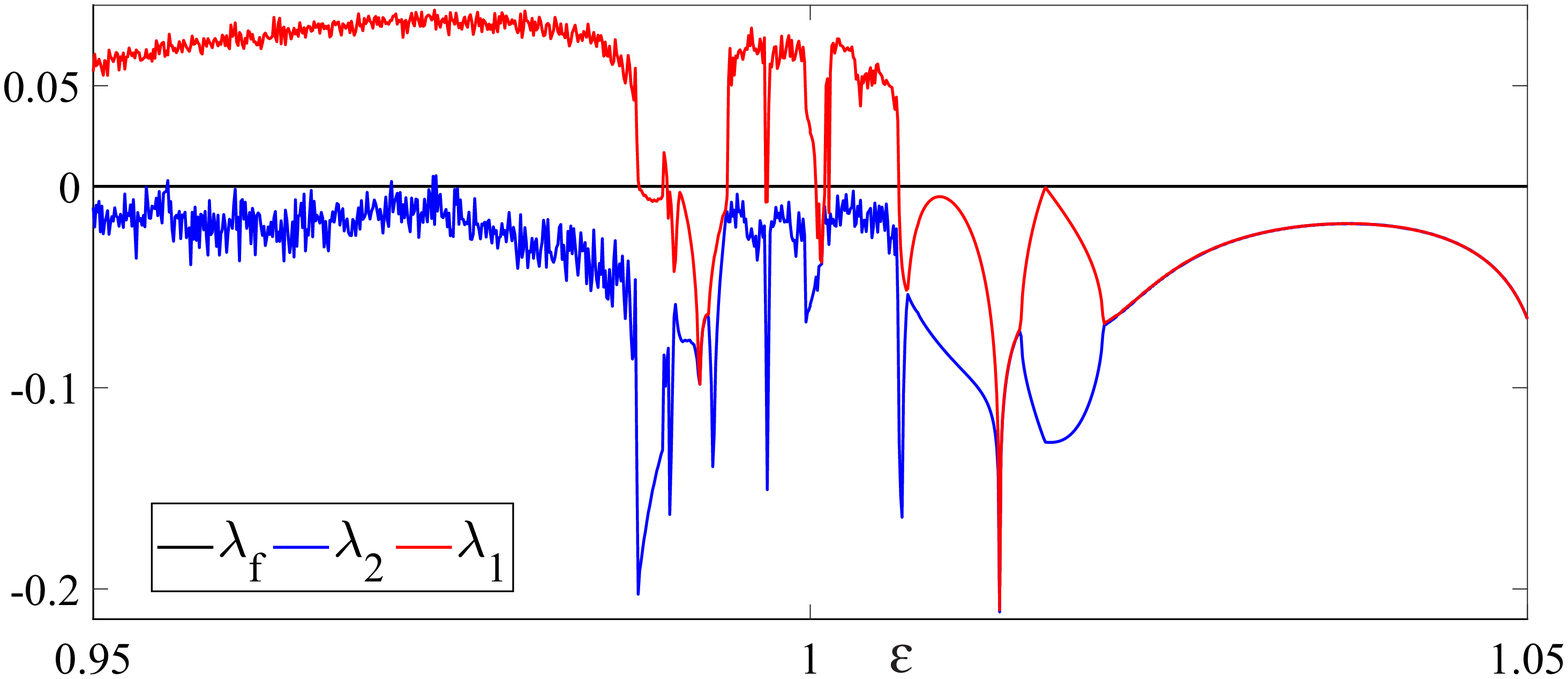}}
\caption{$P_{ac} = 1.35$ MPa, $d/R_{10} = 28$. Dependence of the Lyapunov spectra of attractors on $\varepsilon \in [0.95, 1.05]$.}
\label{fig:map_6}
\end{figure}

Now we compare the results obtained from Fig.-s \ref{fig8}d and \ref{fig9}c: the  synchronous attractors at $d/R_{10} = 11.3$ (point A in Fig. \ref{fig8}a) and $d/R_{10} = 11.47$ (point D in Fig. \ref{fig8}a), continued from the same attractor, evolve in very different ways with respect to changes in $\varepsilon$, depending on the presence or absence of an asynchronous attractor. The monostable synchronous chaotic attractor at $d/R_{10} = 11.3$ is stable to perturbations in $\varepsilon$, while the same synchronous chaotic attractor at  $d/R_{10} = 11.47$ is completely unstable with respect to symmetry breaking. This result confirms our hypothesis that the stability of a synchronous attractor depends mostly on the presence and type of an asynchronous attractors potentially coexisting with the synchronous one, rather than on the type of the synchronous attractor itself.

Finally, we consider the influence of symmetry breaking on another monostable chaotic attractor at  $P_{ac} = 1.35$ MPa, $d/R_{10} = 28$ (substantially weaker coupling). This attractor stays monostable in the entire interval $\varepsilon \in [0.95, 1.05]$, see Fig. \ref{fig:map_6}, according to our numerics. However, this attractor changes its type quite quickly if we change $\varepsilon$. Such a situation can be undesirable for applications, because the type of an actual dynamical regime can differ from the predicted one at symmetrical case. Thus even if an attractor is stable to symmetry breaking, it still does not mean that the corresponding values of the control parameters are necessarily desirable for applications.

\section{Conclusion}
\label{sec:concl}

In this work we have studied the dynamics of two coupled microbubble contrast agents. We have focused on stability of synchronous attractors under the symmetry breaking perturbations. The symmetry breaking has been introduced via variation of parameters corresponding to the ratio of the equilibrium radii of the interacting contrast agents. We have provided several examples of both stable and unstable synchronous attractors.

If a synchronous attractor is monostable, it is very likely to be stable with respect to variations in parameters leading to symmetry breaking. If an asynchronous attractor is a limit cycle, a synchronous attractor coexisting with it, is also likely to be stable. If an asynchronous attractor is quasiperiodic, the synchronous one is also stable with respect to small variations in $\varepsilon$, but exists within a shorter interval of $\varepsilon$, than in the previous case. When an asynchronous attractor is chaotic there are various possibilities. First, a synchronous attractor can be stable with almost the same properties as if it would be coexisting with a quasiperiodic attractor. Second, a synchronous attractor can be unstable with respect to small variations in $\varepsilon$. If an asynchronous attractor is hyperchaotic, the coexisting synchronous attractor is always unstable. Moreover, if a hyperchaotic attractor coexists with another asynchronous attractor, latter one is likely to be unstable under symmetry breaking perturbations as well. On the other hand, hyperchaotic attractors are generally stable with respect to $\varepsilon$ perturbations, which means they are more likely to be observed in physically realistic situation. It can be beneficial for applications, when the desirable dynamical regime is chaotic.

We have demonstrated that the stability of a synchronous attractor is mainly determined by an asynchronous attractor with which it possibly coexists. One of the consequences of this fact is that one should be careful when applying results obtained for synchronous regimes only. For example, in work \cite{Behnia2019} authors carried out an analysis of dynamics in synchronous states for a large ensemble of contrast agents. However, some of found dynamical regimes can be unstable  under small perturbations in the equilibrium radii of bubbles within the ensemble, and, thus, such regimes are not physically realistic. Therefore, we believe that for the studying synchronization of ensembles of contrast agents in particular, and for studying their dynamics in general, it is important to take multistability into account.

Finally, we would like to note, that all observed chaotic and hyperchaotic attractors in the model under consideration are table only from a physical point of view. This means that one can never be sure that the increase in time or accuracy of computation would make the maximal Lyapunov exponent vanish. In the corresponding Lyapunov diagram, one can see, that regions with hyperchaotic behavior alternate with stability windows. The same also concerns other known similar models from various applications, see e.g. Refs. \cite{Stan18, Stan19, Stan20}, demonstrating hyperchaotic attractors. From this point of view, an interesting question arises: which physically relevant systems can have pseudohyperbolic hyperchaotic attractors. \\

This work was supported by Russian Science Foundation grant No. 19-71-10048 (Sec.-s \ref{sec:eqsys}-\ref{sec:one-route}) and by Laboratory of Dynamical Systems and Applications NRU HSE, of the Ministry of science and higher education of the RF grant ag. No. 075-15-2019-1931 (Sec. \ref{sec:2d}).

\end{document}